\documentclass[11pt]{article}
\usepackage{graphicx, colordvi}
\usepackage{amsfonts}
\usepackage{amssymb}
\usepackage{amsthm,cite,color}
\usepackage{dsfont}
\usepackage{epsfig}
\usepackage{mathrsfs}
\usepackage{amsfonts}
\usepackage{amssymb}
\usepackage{amsmath}
\usepackage{amssymb,amsfonts,amsmath,amsthm,cite,color}
\usepackage{dsfont}
\usepackage{epsfig}
\usepackage{mathrsfs}
\usepackage{longtable}
\usepackage{hyperref}
\hypersetup{
colorlinks=true,
linkcolor=cyan,
filecolor=blue,
urlcolor=red,
citecolor=green,
}

\newtheorem{theo}{Theorem}

\newtheorem{lem}[theo]{Lemma}
\newtheorem{exm}[theo]{Example}

\makeatletter \@addtoreset{equation}{section}
\@addtoreset{theo}{section} \makeatother

\newcommand{\bN} { {\mathbb{N}}}

\newcommand{\bZ} { {\mathbb{Z}}}

\newcommand{\bK} { {\mathbb{K}}}

\newcommand{\ord} { {\mathrm{ord\hspace{0.5ex}}}}

\newcommand{\ann} { {\mathrm{ann\hspace{0.5ex}}}}

\newcommand{\Domb}{\mathrm{Domb}}

\def\qed{\hfill \rule{4pt}{7pt}}
\def\pf{\noindent {\it Proof.} }

\begin{document}
\begin{center}

 {\large \bf Rational reductions for holonomic sequences}
\end{center}
\begin{center}
{  Rong-Hua Wang}

  School of Mathematical Sciences\\
   Tiangong University \\
   Tianjin 300387, P.R. China \\
   wangronghua@tiangong.edu.cn
\end{center}

\vskip 6mm \noindent {\bf Abstract.}
Given a holonomic sequence $F(n)$, we characterize rational functions $r(n)$ so that $r(n)F(n)$ can be summable.
We provide upper and lower bounds on the degree of the numerator of $r(k)$ and show the denominator of $r(n)$ can be read from annihilators of $F(k)$.
This illustration provides the so-called rational reductions which can be used to generate new multi-sum equalities and congruences from known ones.

\noindent {\bf Keywords}: rational reduction; holonomic sequence; congruences.

\noindent
\emph{Mathematics Subject Classification}: 11A07, 11B65
\section{Introduction}

In 1990s, Wilf and Zeilberger \cite{WZ1990,Zeilberger1990c} developed the classical WZ theory which can prove and discover new combinatorial identities mechanically.
Zeilberger's algorithm (also known as the method of creative telescoping), designed based on the celebrated Gosper's algorithm \cite{Gosper1978}, is the core of the WZ theory.
The crucial step of algorithm is determining the {\em summability} of a given hypergeometric term $t_n$. That is determining whether $t_n$ can be written as a difference of another hypergeometric term.
During the past three decades, Zeilberger's algorithm has been intensively developed and applied to combinatorics, number theoy and mathematical physics.
At the same time, it is still a attracting subject of ongoing research.
One can consult \cite{CK2017} for a selection of open problems in this context.

The classical Ramanujan-type series for $\pi^{-1}$ are of the form
\begin{equation}\label{eq:r1}
\sum_{n=0}^{\infty}(bn+c)\frac{a(n)}{m^n}
=\frac{\lambda\sqrt{\mu}}{\pi},
\end{equation}
where $b,c,m$ are integers with $bm\neq 0$, $\mu$ is a positive and square free integer, $\lambda$ is a nonzero number, and $a(n)$ is a product of binomial coefficients in $\{
\binom{2n}{n}^3,
\binom{2n}{n}^2\binom{3n}{n},
\binom{2n}{n}^2\binom{4n}{2n},
\binom{2n}{n}\binom{3n}{n}\binom{6n}{3n}
\}.$
Extensive work has been done around applications, $q$-ananlogue forms and $p$-adic analogue forms of Ramaujan-type series. One can consult \cite{BB2010,Berndt1994,CC1988,Guo2020,GL2021,Hamme1997,
HouKraSun2019,Sun2020b,Sun2021,Sun2022} and references therein for more related results.

In 2020, Sun\cite{Sun2020} deduced many new series, similar to the classical Ramanujan-type series, in the form
\begin{equation}\label{eq:rational ramatype}
\sum_{n=0}^{\infty}r(n)
\frac{a(n)}{m^n}
=\frac{\lambda\sqrt{\mu}}{\pi}
\end{equation}
by the method of telescoping.
Here $r(n)$ are rational functions, $a(n)$ and $m$ are give as in \eqref{eq:r1}.
It is straightforward to see that $\frac{a(n)}{m^n}$ is hypergeometric.
Later, Hou and Li \cite{HouLi2021} gave upper and lower bounds on the degree of the numerator of $r(n)$ such that $r(n)t_n$ is summable for a given hypergeometric term $t_n$,
and showed that the denominator of $r(n)$ can be read from the Gosper representation of $t_{n+1}/t_n$. Based on the obtained criteria, they presented a systematic method to construct series like \eqref{eq:rational ramatype}.

There are also many series for $\pi^{-1}$ involving other combinatorial sequences.
For instance, in 2004, Chan, Chan and Liu \cite{CCL2004} and Rogers \cite{Rogers2009} derived
\[
\sum_{n=0}^{\infty}\frac{5n+1}{64^n}\Domb(n)=\frac{8}{\sqrt{3}\pi} \text{ and } \sum_{n=0}^{\infty}(3n+1)\frac{\Domb(n)}{(-32)^n}=
\frac{2}{\pi}
\]
respectively. Here
$
\Domb(n)=\sum_{k=0}^n \binom{n}{k}^2\binom{2k}{k}\binom{2(n-k)}{n-k}
$ is the \emph{Domb numbers}.
Later, Sun presented many interesting conjectures involving the generalized central coefficient $T_n(b,c)$ which denotes the coefficient of $x^n$ in the expansion of $(x^2+bx+c)^n$. For example, Sun \cite{Sun2014,Sun2020b} posed the following conjectural series
\[
\sum_{n=0}^{\infty}\frac{1638n+277}{(-240)^{3n}}\binom{2n}{n}\binom{3n}{n}T_{3n}(62,1)
=\frac{44\sqrt{105}}{\pi},
\]
which was confirmed by J. Wan and W. Zudilin \cite{WanZudilin2012}.
In 2021, Guo and Lian \cite{GL2021} found two series involving $H_n^{(2)}$ corresponding to some classical Ramanujan-type series.
These two identities were proved by Wei \cite{Wei2023}.
Here $H_n^{(m)}=\sum_{k=1}^{n}\frac{1}{k^m}$ is the harmonic number of order $m$.
Recently, more series with summands involving Harmonic numbers has been discussed in \cite{Sun2022,HouSun2023}. Especially, Sun \cite{Sun2022} posed many interesting conjectures on the values of some series involving harmonic numbers of orders not exceeding three.

Let $\bK$ be a field of characteristic $0$. The \emph{annihilator} of a sequence $F(n)$ over $\bK$, denoted by $\ann F(n)$, is defined by
\begin{equation}\label{eq:ann}
\ann F(n):=\left\{L=\sum_{i=0}^{J}a_i(n)\sigma^i\in \bK[n][\sigma]\mid L(F(n))=0\right\},
\end{equation}
where $\sigma$ is the shift operator and $L(F(n))=\sum_{i=0}^{J}a_i(n)F(n+i)$.
A sequence $F(n)$ is called \emph{holonomic} or \emph{P-recursive} if and only if $\ann F(n)\neq \{0\}$.
We call $J$ in \eqref{eq:ann} the order of $L$.
The minimum order of all such recurrences is called the \emph{order} of $F(n)$, denoted by $\ord F(n)$.

By Zeilberger's algorithm, one can check that all afore mentioned sequences such as $\Domb(n)$, $T_n(b,c)$ and $H_n^{(2)}$ are holonomic.
A holonomic sequence $F(n)$ of order $J$ is call \emph{summable} if
\[
F(n)=\Delta\left(\sum_{i=0}^{J-1}u_i(n)F(n+i)\right)
\]
for some rational functions $u_i(n)\in\bK(n)$, $0\leq i\leq J-1$.
Here $\Delta=\sigma-1$ is the difference operator.

It happens frequently that a holonomic sequence $F(n)$ is not summable.
However, when multiplied with a suitable rational function $r(n)$, the product can be summable.
In 1999, Abramov and van Hoeij \cite{AVH1999} provided an algorithm to determine the summaility of $r(n)F(n)$ for any given rational function $r(n)$ and holonomic function $F(n)$.
Suppose $L\in\ann F(n)$ is of order $\ord F(n)$, they proved that $r(n)F(n)$ is summable if and only if the corresponding adjoint equation with respect to $L$ has rational solutions.


In this paper, we will characterise rational functions $r(n)=a(n)/b(n)$ such that $r(n)F(n)$ is summable for any given holonomic sequence $F(n)$.
Here $a(n),b(n)$ are polynomials over $\bK$ and $\gcd (a(n),b(n))=1$.
More precisely, we first provide an upper bound and a lower bound on the degree of $p(n)\in\bK[n]$ if $p(n)F(n)$ is summable.
The upper bound given here can be better than the one provided in \cite{HouLi2021} for some hypergeometric terms.
Then we will show the denominator $b(n)$ can be read from annihilators of $F(n)$.
By choosing appropriate $b(n)$, one can generate new multi-sum identities and congruences mechanically from known ones.

The paper is organized as follows.
Section 2 is devoted to providing upper and lower bounds on the degree of $p(n)\in\bK[n]$ such that $p(n)F(n)$ is summable. In Section 3, we first characterize the denominator of rational functions $r(n)$ with $r(n)F(n)$ being summable, and then introduce two kinds of rational reductions. Examples of discovery and proof of new equalities and congruences automatically are provided in Section 4.

\section{Summablity of $p(n)F(n)$}
Let $F(n)$ be a holonomic sequence, in this section, we will provide upper and lower bounds on the degree of $p(n)\in\bK[n]$ such that $p(n)F(n)$ is summable.

For any recurrence operator $L=\sum_{i=0}^{J}a_i(n)\sigma^i$, the \emph{adjoint} of $L$ is defined by
$L^{\ast}=\sum_{i=0}^{J}\sigma^{-i}a_i(n)$, and thus
\[
L^{\ast}(x(n))=\sum_{i=0}^{J}a_i(n-i)x(n-i)
\]
for any $x(n)\in\bK[n].$
When $L\in\ann F(n)$, by \cite[eq (1.6)]{CHJ2011} or \cite[Proposition 3.2]{Hoeven2018} we know
\begin{equation}\label{eq: summable}
L^{\ast}(x(n))F(n)=\Delta
\left(-\sum_{i=0}^{J-1}u_i(n)F(n+i)\right),
\end{equation}
where
\begin{equation}\label{eq:u}
u_i(n)=\sum_{j=1}^{J-i}a_{i+j}(n-j)x(n-j).
\end{equation}
That means $L^{\ast}(x(n))F(n)$ is summable for any $x(n)\in\bK[n]$.
From equality \eqref{eq: summable} we obtain
\begin{equation}\label{eq:rec ann}
\sum_{n=\ell}^{k-1}L^{\ast}(x(n)) F(n)
=\left(\sum_{i=0}^{J-1}u_i(\ell)F(\ell+i)\right)
-\left(\sum_{i=0}^{J-1}u_i(k)F(k+i)\right).
\end{equation}

Given $L=\sum\limits_{i=0}^{J}a_i(n)\sigma^i\in\ann F(n)$.
\emph{The degree} of $L$, written as $\deg L$, is given by
\begin{equation}\label{eq:d and bk}
\deg L=\max_{0\leq k \leq J} \{\deg b_k(n)-k\},
\end{equation}
where
\begin{equation}\label{eq:d and bk1}
b_k(n)=\sum_{j=k}^{J}\binom{j}{k}a_{J-j}(n+j-J),\quad 0\leq k\leq J.
\end{equation}
Let $d=\deg L$ and
\begin{equation}\label{eq:f(s)}
f(s)=\sum_{k=0}^{J}[n^{d+k}](b_k(n))s^{\underline{k}},
\end{equation}
where
$[n^{d+k}](b_k(n))$ denotes the coefficient of $n^{d+k}$ in $b_k(n)$ and $s^{\underline{k}}$ denotes the falling factorial $s(s-1)\cdots (s-k+1)$.
Let
\begin{equation}\label{eq:nonnegative roots}
R_{L}=\{s\in\bN \mid f(s)=0\}.
\end{equation}
Then $L$ is called \emph{degenerated} if $R_{L}\neq \emptyset$ and \emph{nondegenerated} otherwise.

In 2023, Wang and Zhong \cite[Lemma 2.3]{WZ2023AAM} proved that for any nonzero polynomial $x(n)\in\bK[n]$, we have
\begin{equation}\label{eq:deg L(p(n))}
\deg L^{\ast}(x(n))\left\{
     \begin{array}{ll}
       < \deg L+\deg x(n), & \hbox{if $L$ is degenerated and $\deg x(n)\in R_{L}$,} \\
       =\deg L+\deg x(n), & \hbox{otherwise.}
     \end{array}
   \right.
\end{equation}

The following theorem shows that an upper bound on the degree of $p(n)$ such that $p(n)F(n)$ is summable can be given by the degree and continued zero index of $L$.
The \emph{continued zero index} of $L$, denoted by $C_L$, is defined as follows:
if $L^{\ast}(1)\neq 0$, then $C_L=0$; if $L^{\ast}(1)=0$, then $C_L$ is the positive integer such that $L^{\ast}(n^{C_L})\neq 0$ and
$L^{\ast}(n^{i})=0$ if $0\leq i \leq C_L-1$.

\begin{theo}\label{th:upper bound}
Suppose $F(n)$ is a holonomic sequence and $L=\sum\limits_{i=0}^{J}a_i(n)\sigma^i\in\ann F(n)$.
Then there exists a nonzero polynomial $p(n)$ such that $p(n)F(n)$ is summable and
\begin{equation}\label{eq:upper bound}
\deg p(n)\leq d+C_L,
\end{equation}
where $d=\deg L$ and $C_L$ is the continued zero index of $L$.
\end{theo}
\pf By the definition of $C_L$, we have $L^{\ast}(n^{C_L})\neq 0$.
With the help of \eqref{eq:deg L(p(n))} and \eqref{eq: summable}, one can derive that
\begin{equation}\label{eq:L(x(n))}
\deg L^{\ast}(n^{C_L})\leq d+C_L
\end{equation}
and $L^{\ast}(n^{C_L})F(n)$ is summable.
\qed

The inequality \eqref{eq:upper bound} provides an upper bound.
It may happen that there exists a polynomial $p(n)$ such that $p(n)F(n)$ is summable and $\deg p(n)<d+C_L$.
\begin{exm}
Let
\[
F(n)=\frac{(-1)^n f_n}{n(n-1)},
\]
where $f_n$ denotes the Franel number $f_n=\sum_{k=0}^{n}\binom{n}{k}^3$.
By Zeilberger's algorithm, one can see
\[
L=8n(n^2-1)-n(7n^2+21n+16)\sigma-(n+2)^3\sigma^2
\in\ann F(n).
\]
Then $J=\ord F(n)=2$, $d=\deg L=2$, $C_L=0$, and thus the upper bound given in
\eqref{eq:upper bound} is $2$.
Take $p(n)=2-3n$, one can check that
\[
(2-3n)F(n)=\Delta\left(  (8n^3-5n-2)F(n)+(n+1)^3 F(n+1)\right).
\]
Then we know $(2-3n)F(n)$ is summable while $\deg (2-3n)=1<d$.
\end{exm}

The following example shows that the continued zero index $C_L$ in \eqref{eq:upper bound} is necessary.
\begin{exm}\label{ex:harmonic}
Let $H_n=\sum_{k=1}^{n}\frac{1}{k}$ denote the Harmonic number and
\[
F(n)=\frac{H_n}{n(n+1)}.
\]
By Zeilberger's algorithm, one can see
\[
L=n(n+1)^2-(n+1)(n+2)(2n+3)\sigma+(n+2)^2(n+3)\sigma^2
\in\ann F(n).
\]
Then $J=\ord F(n)=2$, $d=\deg L=1$ and $C_L=1$.
Let $a_0(n)=n(n+1)^2$ and $a_2(n)=(n+2)^2(n+3)$.
It is easy to check that $\gcd(a_0(n),a_2(n+h))=1$ for any $h\in\bN=\{0,1,2,\ldots\}$.
Note that
\[
(n^2 + n)F(n)=\Delta\left(  n (n^3 + 2 n^2 -1 )F(n)-(n-1) ( n+1)^2 (n+2) F(n+1)\right).
\]
Then we know $(n^2 + n)F(n)$ is summable, here $\deg (n^2 + n)=d+C_L$=2.
Since $C_L=1$, Theorem 2.4 and Lemma 2.5 in \cite{WZ2023AAM} lead to the fact that there exists no polynomial $p(n)$ with $\deg p(n)\leq 1$ such that $p(n)F(n)$ is summable.
\end{exm}
Example \ref{ex:harmonic} also illustrates that for some holonomic sequences, the upper bound given in Theorem \ref{th:upper bound} can also be a lower bound.
Next we show this happens frequently.
To this aim, we first provide two lemmas on $C_L$.
\begin{lem}\label{lem:C_L<J}
Let $L=\sum\limits_{i=0}^{J}a_i(n)\sigma^i$. Then the continued zero index $C_L\leq J$.
\end{lem}
\pf
Note that $f(s)$ given by  \eqref{eq:f(s)} is a polynomial in $s$ and $\deg f(s)\leq J$.
Thus the set $R_{L}$ defined in \eqref{eq:nonnegative roots} contains at most $J$ elements. The proof then follows directly from \eqref{eq:deg L(p(n))} and the definition of $C_L$.
\qed
\begin{lem}\label{lem:C_L=0}
If $L=\sum\limits_{i=0}^{J}a_i(n)\sigma^i$ is nondegenerated. Then $C_L=0$.
\end{lem}
\pf
As $L$ is nondegenerated, we know $\deg L^{\ast}(x(n))=d+\deg x(n)$ by \eqref{eq:deg L(p(n))}.
Here $d=\deg L$.
Then $\deg L^{\ast}(n^{C_L})=d+C_L\geq 0$ as $L^{\ast}(n^{C_L})\neq 0$, and thus
$d \geq -C_L\geq -J$ by Lemma \ref{lem:C_L<J}.
Note that $\deg L^{\ast}(1)=d\geq -J$, one can see $ L^{\ast}(1)\neq 0$.
Hence $C_L=0$.
\qed
\begin{theo}
Let $L=\sum\limits_{i=0}^{J}a_i(n)\sigma^i\in\ann F(n)$ with  $a_0(n)a_J(n)\neq 0$ and $J>0$ be the order of $F(n)$.
Suppose $\gcd(a_0(n),a_J(n+h))=1$ for any $h\in\bN$ and $L$ is nondegenerated.
If $p(n)$ is a nonzero polynomial such that $p(n)F(n)$ is summable, then
\[
\deg p(n)\geq \deg L+C_L.
\]
\end{theo}
\pf
Since $\gcd(a_0(n),a_J(n+h))=1$ for any $h\in\bN$, we have $p(n)F(n)$ is summable if and only if $p(n)=L^{\ast}(x(n))$ for some $x(n)\in\bK[n]$ by \cite[Theorem 2.4]{WZ2023AAM}.
As $L$ is nondegenerated, we have
\[
\deg p(n)=\deg L^{\ast}(x(n))=\deg L+\deg x(n)\geq \deg L+C_L
\]
by Lemma \ref{lem:C_L=0}
\qed

Since hypergeometric terms are all holonomic, Theorem \ref{th:upper bound} also provides an upper bound for hypergeometric terms. In general, the upper bound provided in \cite{HouLi2021} may be superior than \eqref{eq:upper bound} since Q.-H. Hou and G.-J. Li \cite{HouLi2021} considered the Gosper representation of $F(n)$.
The following example shows when $F(n)$ contains no polynomial part, the upper bound given by Theorem \ref{th:upper bound} can be better.
\begin{exm}
Let
\[
F(n)=\frac{\binom{2n}{n}^4}{(2n-1)^4 256^n}.
\]
It is easy to check that
\[
L=(2n-1)^4-16(n+1)^4\sigma
\in\ann F(n).
\]
Then $J=\ord F(n)=1$, $d=\deg L=3$ and $C_L=0$.
Hence the upper bound given in \eqref{eq:upper bound} is $d=3$ while the one presented in \cite{HouLi2021} is $4$.
\end{exm}
\section{Summablity of $r(n)F(n)$}
In this section, we will characterise the denominator of rational function $r(n)$ such that $r(n)F(n)$ is summable.
The criteria motivate the discovery of the so-called \emph{rational reduction}, which can be used to generated new equalities from known ones.
\begin{theo}\label{th:rational summable}
Let $L=\sum\limits_{i=0}^{J}a_i(n)\sigma^i\in\ann F(n)$ with $a_0(n)a_J(n)\neq 0$ and $J>0$ be the order of $F(n)$.
Let $a(n),b(n)\in\bK[n]$ such that $\frac{a(n)}{b(n)}F(n)$ is summable.
Suppose that for $\forall h\in\bN$, we have
\begin{equation}\label{eq:a0J}
\gcd(a_0(n),a_J(n+h))=\gcd(b(n),b(n+J+h))=1
\end{equation}  and
\begin{equation}\label{eq:condition}
\gcd(a_0(n),b(n+J+h))=\gcd(b(n),a_J(n+h))=1.
\end{equation}
Then $b(n)\mid a(n)$.
\end{theo}
\pf
Let $G(n)=\frac{F(n)}{b(n)}$. Then $a(n)G(n)$ is summable.
Since $L=\sum\limits_{i=0}^{J}a_i(n)\sigma^i\in\ann F(n)$, one can check that
\[
L_1=\sum\limits_{i=0}^{J}a_i(n)b(n+i)\sigma^i\in\ann G(n).
\]
Take $b_i(n)=a_i(n)b(n+i)$, $0\leq i \leq J$.
When conditions in Theorem \ref{th:rational summable} are satisfied, we have
$\gcd(b_0(n),b_J(n+h))=1$ for any nonnegative integer $h$.
Thus by \cite[Theorem2.4]{WZ2023AAM} we derive that
\[
a(n)=L_1^{\ast}(x(n))=b(n)\sum\limits_{i=0}^{J}a_i(n-i)x(n-i)=b(n)L^{\ast}(x(n)),
\]
for some $x(n)\in\bK[n]$.
This completes the proof.
\qed

By Theorem \ref{th:rational summable}, if $\frac{a(n)}{b(n)}F(n)$ is summable with $b(n)\nmid a(n)$ and conditions \eqref{eq:a0J} holds, then $b(n)$ must contain factors from $a_0(n-J-h)$ or $a_J(n+h)$ for some nonnegative integer $h$.

To proceed, some notations are needed.
Suppose $A(n)\in\bK[n]$ and $I\in\bZ$.
The \emph{shift product operator of order $I$}, denoted by $SP_I$, is defined by
\begin{equation}\label{eq:SP}
SP_{I}(A(n))=\left\{
     \begin{array}{lll}
       \prod\limits_{j=1}^{I}A(n+j), & \hbox{when $I>0$,} \\[10pt]
       1, & \hbox{when $I=0$,} \\[5pt]
       \prod\limits_{j=1}^{-I}A(n-j), & \hbox{when $I<0$.}
     \end{array}
   \right.
\end{equation}
Let $L=\sum\limits_{i=0}^{J}a_i(n)\sigma^i$, $a_i(n)\in\bK[n]$. Define
\begin{equation}\label{eq:d_L}
d_L=\max_{0\leq i \leq J}\{\deg a_i(n)\}.
\end{equation}
The following lemma shows the connection between $d_L$ and $\deg L$ is closely related to the degenerality of $L$.
\begin{lem}
Given $L=\sum\limits_{i=0}^{J}a_i(n)\sigma^i$, $a_i(n)\in\bK[n]$, and $d_L$ be defined by
\eqref{eq:d_L}.
Then $\deg L\leq d_L$, and $L$ is nondegenerated if $\deg L=d_L$.
\end{lem}
\pf By the definition of $b_k(n)$ in \eqref{eq:d and bk1}, it is easy to see $\deg b_k(n)\leq d_L$.
Then
$$\deg L=\max_{0\leq k\leq J}\{\deg b_k(n)-k\}\leq \max_{0\leq k\leq J}\{d_L-k\}=d_L.$$
When $\deg L=d_L$, we know $\deg b_0(n)=d_L$.
thus
\begin{equation*}
f(s)=\sum_{k=0}^{J}[n^{\deg L+k}](b_k(n))s^{\underline{k}}=[n^{d_L}](b_0(n))
\end{equation*}
is a nonzero constant.
So $L$ is nondegenerated.
\qed

When $\deg L=d_L$, the operator $L$ is called \emph {strongly nondegenerated}, as the condition $\deg L=d_L$ is not necessary for $L$ to be nondegenerated.
\begin{exm}\label{exm:Domb/16^n}
Let $F(n)=\Domb(n)/16^n$, Zeilberger's algorithm leads to
\[
L=8 (2 + n)^3 \sigma^2 - (3 + 2 n) (12 + 15 n + 5 n^2) \sigma + 2 (1 + n)^3\in\ann F(n).
\]
One can check that
$\deg L=2$ and $L$ is nondegenerated while $d_L=3>\deg L$.
\end{exm}
Assume $L=\sum\limits_{i=0}^{J}a_i(n)\sigma^i\in \ann F(n)$ is nondegenerated.
For any polynomial $p(n)\in\bK[n]$, by equality \eqref{eq:deg L(p(n))} we can decompose $p(n)$ as
\[
p(n)=q(n)+\tilde{p}(n),
\]
such that $q(n),\tilde{p}(n)\in\bK[n]$, $q(n)F(n)$ is summable and $\deg \tilde{p}(n)<\deg L$, utilizing the division algorithm.
The above decomposition is called {\em polynomial reduction} of $p(n)$ with $L$ in \cite{WZ2023AAM}.

Next, we will introduce two kinds of rational reductions which can be used after polynomial reduction.
\begin{theo}\label{th:main1}
Let $F(n)$ be a holonomic sequence, $L=\sum\limits_{i=0}^{J}a_i(n)\sigma^i\in\ann F(n)$ is strongly nondegenerated and
$
a_0(n)=A_0(n)\tilde{A}_0(n),
$
where $A_0(n),\tilde{A}_0(n)\in\bK[n]$.
For any $p(n)\in\bK[n]$, if positive integer $I\geq J$, we can find a polynomial $\tilde{p}(n)$ and a holonomic sequence $T(n)$ such that
\begin{equation}\label{eq:rational decom1}
p(n)F(n)=\frac{\tilde{p}(n)}
{SP_{-I}(A_0(n)))}F(n)
+\Delta(T(n)),
\end{equation}
and $\deg \tilde{p}(n)<\deg L+(J-1)\deg A_0(n)$.
\end{theo}
\pf
Take $G(n)=\frac{F(n)}{SP_{-I}(A_0(n))}$.
Since $L=\sum\limits_{i=0}^{J}a_i(n)\sigma^i\in\ann F(n)$, it is straightforward to check that
\[
\sum_{i=0}^{J}a_i(n)SP_{-I}(A_0(n+i))\sigma^i\in\ann G(n).
\]
Note that
\begin{align*}
&\sum_{i=0}^{J}a_i(n)SP_{-I}(A_0(n+i))\sigma^i \\
=& \tilde{A}_0(n)\prod_{j=0}^{I-J}A_0(n-j)
  \prod_{j=I-J+1}^{I}A_0(n-j)\\
& + \sum_{i=1}^{J}a_i(n)\prod_{j=1}^{i-1}A_0(n+j)
  \prod_{j=0}^{I-J}A_0(n-j)
  \prod_{j=I-J+1}^{I-i}A_0(n-j)\sigma^i.
\end{align*}
Then we have
\[L_1=\tilde{A}_0(n)\prod_{j=I-J+1}^{I}A_0(n-j)
  + \sum_{i=1}^{J}a_i(n)\prod_{j=1}^{i-1}A_0(n+j)
  \prod_{j=I-J+1}^{I-i}A_0(n-j)\sigma^i
\]
is an annihilator of $G(n)$.
Note that $L_1$ is strongly nondegenerated as $L$ is strongly nondegenerated and the shift opeartor $\sigma^j$ does not change the leading coefficient. Hence
\[
\deg L_1=\deg L+(J-1)\deg A_0(n).
\]
The polynomial reduction with $L_1$ leads to the decompostition
\[
p(n)SP_{-I}(A_0(n))=\tilde{p}(n)+q(n),
\]
such that $\tilde{p}(n)<\deg L_1$ and $q(n)G(n)=\Delta(T(n))$ for some holonomic sequence $T(n)$.
Then
\begin{align*}
p(n)F(n)
=& p(n)SP_{-I}(A_0(n))G(n)=(\tilde{p}(n)+q(n))G(n) \\
=& \frac{\tilde{p}(n)}{SP_{-I}(A_0(n))}F(n)+\Delta(T(n)).
\end{align*}
This completes the proof.
\qed
\begin{theo}\label{th:main2}
Let $F(n)$ be a holonomic sequence, $L=\sum\limits_{i=0}^{J}a_i(n)\sigma^i\in\ann F(n)$ is strongly nondegenerated and
$
a_J(n)=A_J(n)\tilde{A}_J(n),
$
where $A_J(n),\tilde{A}_J(n)\in\bK[n]$.
For any $p(n)\in\bK[n]$, if positive integer $I\geq J$, we can find a polynomial $\tilde{p}(n)$ and a holonomic sequence $T(n)$ such that
\begin{equation}\label{eq:rational decom2}
p(n)F(n)=\frac{\tilde{p}(n)}
{SP_{I}(A_J(n-J)))}F(n)
+\Delta(T(n)),
\end{equation}
and $\deg \tilde{p}(n)<\deg L+(J-1)\deg A_J(n)$.
\end{theo}
\pf
Take $G(n)=\frac{F(n)}{SP_{I}(A_J(n-J))}$.
Note that
\[
\sum_{i=0}^{J}a_i(n)SP_{I}(A_J(n-J+i))\sigma^i\in\ann G(n).
\]
and  that
\begin{align*}
&\sum_{i=0}^{J}a_i(n)SP_{I}(A_J(n-J+i))\sigma^i \\
&= \sum_{i=0}^{J-1}a_i(n)\prod_{j=1}^{J-i-1}A_J(n-j)
  \prod_{j=0}^{I-J}A_J(n+j)
  \prod_{j=I-J+1}^{I-J+i}A_J(n+j)\sigma^i\\
&+\tilde{A}_J(n)\prod_{j=0}^{I-J}A_J(n+j)
  \prod_{j=I-J+1}^{I}A_J(n+j)\sigma^J.
\end{align*}
Thus one can see
\[\sum_{i=0}^{J-1}a_i(n)\prod_{j=1}^{J-i-1}A_J(n-j)
   \prod_{j=I-J+1}^{I-J+i}A_J(n+j)\sigma^i
+\tilde{A}_J(n)
  \prod_{j=I-J+1}^{I}A_J(n+j)\sigma^J
  \]
falls into $\ann G(n)$.
The proof then is similar to that of Theorem \ref{th:main1}.\qed

The decompositions \eqref{eq:rational decom1} and \eqref{eq:rational decom2} will be called the \emph{ rational reduction} to $p(n)$ with $L$.
In Theorem \ref{th:main1} and Theorem \ref{th:main2}, the operator $L$ is assumed to be strongly nondegenerated.
When this is not the case, all discussions in the proof go well except the determination of $\deg L_1$.
Since $\deg L_1$ is free of $I$ and
$\deg(p(n)SP_{I}(A(n)))=\deg p(n)+I \cdot \deg A(n)$ for any polynomial $A(n)$.
When $I$ is large enough, one can still obtain equalities as \eqref{eq:rational decom1} or \eqref{eq:rational decom2}.
Example \ref{ex2} in the next subsection is a concrete illustration.

\section{Applications of rational reductions}
In this section,  we will show how rational reduction provides an algorithmic way to generate new equalities and congruences from known ones.
We will take Domb numbers as examples to illustrate the process.

In 2009, Rogers \cite{Rogers2009} derived
\begin{equation}\label{eq:Domb -32}
\sum_{n=0}^{\infty}(3n+1)\frac{\Domb(n)}{(-32)^n}=
\frac{2}{\pi}.
\end{equation}
Sun \cite{Sun2021} proposed series of conjectured identities with the summand in the form of $p(n)\frac{\Domb(n)}{(-32)^n}$, where $p(n)$ is a polynomial and $\deg p(n)\geq 3$.
The conjectures are confirmed recently in \cite{WZ2023AAM} by the polynomial reduction.
Next, we provide several equalities involving $r(n)\frac{\Domb(n)}{(-32)^n}$ with $r(n)$ being rational using the rational reductions generated in Section 3.
\begin{theo}
The following equalities holds:
\begin{align}
 & \sum_{n=0}^{\infty}\frac{27 + 103 n + 141 n^2 + 78 n^3 + 15 n^4}{(n+1)^2(n+2)^2}\frac{\Domb(n)}{(-32)^n}=
\frac{18}{\pi}\label{eq:1} \\
 &
\sum_{n=0}^{\infty}\frac{239 + 807 n + 993 n^2 + 582 n^3 + 165 n^4 + 18 n^5}{(n+1)^3(n+2)^3}\frac{\Domb(n)}{(-32)^n}=
80 - \frac{162}{\pi}\label{eq:2}  \\
& \sum_{n=2}^{\infty}\frac{2 + 5 n - 9 n^2 - 21 n^3 + 39 n^4}{n^2(n-1)^2}\frac{\Domb(n)}{(-32)^n}=
-9\left(\frac{2}{\pi}-\frac{11}{12}\right) \label{eq:3} \\
&
\sum_{n=2}^{\infty}\frac{-26 - 21 n + 105 n^2 - 9 n^3 - 147 n^4 + 306 n^5}{n^3(n-1)^3}\frac{\Domb(n)}{(-32)^n}=
\frac{162}{\pi}-\frac{217}{8}\label{eq:4} \\
&
\sum_{n=0}^{\infty}\frac{5729 + 8701 n + 5895 n^2 + 1879 n^3 + 228 n^4}{(n+1)^2(n+2)^2(n+3)^2}\frac{\Domb(n)}{(-32)^n}=
\frac{486}{\pi}\label{eq:5} .
\end{align}
\end{theo}
\proof
Let $F(n)=\Domb(n)/(-32)^n$.
By Zeilberger's algorithm, we find
$$L=
(n+1)^3+(2n+3) (5 n^2+15n+12)\sigma+16 (n+2)^3\sigma^2\in\ann F(n).
$$

If we take $A_2(n)=(n+2)^2$ which divides $(n+2)^3$ and $I=2$.
Then
\[
G(n)=\frac{F(n)}{SP_{2}(A_2(n-2))}=\frac{F(n)}{(n+1)^2(n+2)^2}.
\]
By Zeilberger's algorithm, we find that $L_1=\sum_{i=0}^{2} a_i(n) \in\ann G(n)$ with
$a_2(n)=16 (2 + n) (3 + n)^2 (4 + n)^2$, $a_1(n)=(3 + n)^2 (3 + 2 n) (12 + 15 n + 5 n^2)$, and $a_0(n)=(1 + n)^5.$
One can check that $L_1$ is strongly nondegenerated and $\deg L_1=5$.
So $P_s(n)=L_1^{\ast}(n^s)$ is polynomial and $\deg P_s(n)=s+5$ for any $s\ge 0$.
Note that $(3n+1)F(n)$=$q(n)G(n)$, where $q(n)=(3n+1)(n+1)^2(n+2)^2$.
Since $\deg q(n)=5=\deg L_1$, by polynomial reduction to $q(n)$ with $L_1$, we have
\begin{equation}\label{eq:-32-1}
q(n)=\frac{1}{9} (27 + 103 n + 141 n^2 + 78 n^3 + 15 n^4)-\frac{1}{9}P_0(n).
\end{equation}
From equality \eqref{eq: summable}, we know
\[
P_0(n)G(n)=\Delta
\left(-u_0(n)G(n)-u_1(n)G(n+1)\right)
\]
with $u_i(n)=\sum_{j=1}^{2-i}a_{i+j}(n-j)$.
Multiplying both sides of \eqref{eq:-32-1} with $G(n)$ and then summing over $n$ from $0$ to $\infty$ derives equality \eqref{eq:1} from \eqref{eq:Domb -32}.

If we take $A_2(n)=(n+2)^3$ and $I=2$, similar discussions lead to equality
\eqref{eq:2}.
Taking $A_0(n)=(n+1)^2$ or $A_0(n)=(n+1)^3$ and $I=2$ will arrive at equalities
\eqref{eq:3}
and
\eqref{eq:4}
respectively. Note that the summations have to start from $2$.
Equality \eqref{eq:5} follows by taking $A_2(n)=(n+2)^2$ and $I=3$.
We can obtain infinitely many new series by increasing $I$.
\qed

In 2013, Sun \cite{Sun2013} conjectured that when $p$ is a prime and $p\equiv 1\pmod 3$,
\begin{equation}\label{DombCong}
\sum_{n=0}^{p-1}(3n+1)\frac{\Domb(n)}{16^n}\equiv 0 \pmod {p^2}.
\end{equation}
This conjecture was confirmed and generalized by Mao and Liu \cite{MaoLiu2022} recently.
By the rational reduction and congruence \eqref{DombCong}, we can discover new arithmetic properties of Domb numbers.
\begin{theo}\label{ex2}
Let $p$ be a prime and $p\equiv 1\pmod 3$. Then
\begin{equation}\label{DombCong2}
\sum_{n=2}^{p-1}\frac{(n+1)^2}{n(n-1)}\cdot\frac{\Domb(n)}{16^n}\equiv \frac{3}{2} \pmod {p^2}.
\end{equation}
\end{theo}
\pf Let $F(n)=\Domb(n)/16^n$, we know by Example \ref{exm:Domb/16^n} that
\[
L=8 (2 + n)^3 \sigma^2 - (3 + 2 n) (12 + 15 n + 5 n^2) \sigma + 2 (1 + n)^3\in\ann F(n)
\]
and that $L$ is not strongly nondegenerated.
Taking $A_0(n)=n+1$ and $I=2$, we have $G(n)=\frac{F(n)}{SP_{-I}(A_0(n))}=\frac{F(n)}{n(n-1)}$ and
\[
L_1=8 (2 + n)^4 \sigma^2 - n (3 + 2 n) (12 + 15 n + 5 n^2) \sigma +
 2 (-1 + n) n (1 + n)^2\in\ann G(n).
\]
One can check that $L_1$ is nondegenerated and $\deg L_1=3$.
By polynomial reduction with respect to $L_1$, we have
\begin{equation}\label{eq:n(n-1)(3n+1)}
n(n-1)(3n+1)=-L_1^{\ast}(n^0)+2(n+1)^2.
\end{equation}
Equality \eqref{eq:rec ann} leads to
\begin{align*}
\sum_{n=2}^{p-1}L_1^{\ast}(n^0) G(n)
=\ &\sum_{i=0}^{1}u_i(2)G(2+i)-\sum_{i=0}^{1}u_i(p)G(k+p)\\
=\ & 5-2 p^2 (2 - 3 p + p^2) G(-1 + p) + 8 p^4 G(p) \\
\equiv\ &5\pmod {p^2},
\end{align*}
where $u_0(n)=2 + 7 n + 6 n^2 - 5 n^3 - 2 n^4$ and
$u_1(n)=8 (1 + n)^4$.
Multiplying both sides of \eqref{eq:n(n-1)(3n+1)} with $G(n)$ and then summing over $n$ from $2$ to $p-1$ arrives at
\[
\sum_{n=2}^{p-1}\frac{(n+1)^2}{n(n-1)}\cdot \frac{\Domb(n)}{16^n}
\equiv\frac{1}{2}\left( 5+\sum_{n=2}^{p-1}(3n+1)\frac{\Domb(n)}{16^n}\right)\equiv
\frac{3}{2} \pmod {p^2},
\]
when $p\equiv 1 \pmod 3$.
\qed

\noindent \textbf{Acknowledgments.}
This work was supported by the National Natural Science Foundation of China  (No. 12101449, 12271511) and the Natural Science Foundation of Tianjin, China (No. 22JCQNJC00440).


\begin{thebibliography}{10}
\bibitem{AVH1999}
S.A. Abramov and M. van Hoeij.
\newblock{Integration of solutions of linear functional equations.}
\newblock{\em Integral Transform. Spec. Funct.}, 8(1–2)(1999): 3--12.

\bibitem{BB2010}
N.D. Baruah and B.C. Berndt.
\newblock Eisenstein series and Ramanujan-type series for $1/\pi$.
\newblock{\em Ramanujan J.}, 23(2010),17--44.

\bibitem{Berndt1994}
 B.C. Berndt.
\newblock Ramanujan’s Notebooks.
\newblock {\em Part IV}, Springer, NewYork, 1994.
\bibitem{CCL2004}
H.H. Chan, S.H. Chan and Z. Liu.
\newblock Domb's numbers and Ramanujan-Sato type series for $1/\pi$.
\newblock {\em Adv. Math.}, 186(2004), 396--410.
\bibitem{CK2017}
S. Chen and M. Kauers.
\newblock Some open problems related to creative telescoping.
\newblock{\em J. Syst. Sci. Complex}, 30(1)(2017): 154--172.
\bibitem{CHJ2011}
W.Y.C. Chen, Q.-H. Hou and H.-T. Jin.
\newblock The Abel-Zeilberger algorithm.
\newblock {\em Electron. J. Combin.}, 18(2)(2011), P17.
\bibitem{CC1988}
D.V. Chudnovsky and G.V. Chudnovsky.
\newblock Approximations and complex multiplication according to Ramanujan.
In \newblock{ Ramanujan Revisited}, Academic Press, Boston, MA, 1988, 375--472.
\bibitem{Gosper1978}
R.W. Gosper.
\newblock Decision procedure for indefinite hypergeometric summation.
\newblock {\em Proc. Natl. Acad. Sci. U.S.A.}, 75(1):40--42, 1978.
\bibitem{Guo2020}
V.J.W. Guo.
\newblock $q$-Analogues of three Ramanujan-type formulas for $1/\pi$.
\newblock {\em Ramanujan J.}, 52(2020), 123--132.
\bibitem{GL2021}
V.J.W. Guo and X. Lian.
\newblock Some $q$-congruenceson double basic hypergeometric sums.
\newblock {\em J. Difference Equ. Appl.}, 27 (2021), 453--461.


\bibitem{HouKraSun2019}
Q.-H. Hou, C. Krattenthaler and Z.-W. Sun.
\newblock On $q$-analogues of some series for $\pi$ and $\pi ^2$.
\newblock {\em Proc. Amer. Math. Soc.}, 5(2019), 1953--1961.

\bibitem{HouLi2021}
Q.-H. Hou and G.-J. Li.
\newblock Gosper summability of rational multiples of hypergeometric terms.
\newblock {\em J. Difference Equ. Appl.}, 27(2021), 1723--1733.
\bibitem{HouSun2023}
Q.-H. Hou and Z.-W. Sun.
\newblock Taylor coefficients and series involving  harmonic numbers.
\newblock {arXiv:2310.03699v2}.
\bibitem{MaoLiu2022}
G.-S. Mao and Y. Liu.
\newblock Proof of some conjectural congruences involving Domb numbers and binary quadratic forms.
\newblock {\em J. Math. Anal. Appl.},  516 (2022), 126493.

\bibitem{Rogers2009}
M. D. Rogers.
\newblock New $_5 F_4$ hypergeometric transformations, three-variable Mahler measures, and formulas for $1/\pi$.
\newblock {\em Ramanujan J.}, 18 (2009), 327-340.
\bibitem{Sun2013}
Z.-W. Sun.
\newblock  Number Theory and Related Area.
\newblock {\em (eds., Y. Ouyang, C. Xing, F. Xu and P. Zhang), Adv. Lecr. Math.} 27, Higher Education Press and International Press, Beijing-Boston, 2013, pp. 149–197
\bibitem{Sun2014}
Z.-W. Sun.
\newblock On sums related to central binomial and trinomial coefficients.
\newblock in: {\em  Combinatorial and Additive Number Theory: CANT 2011 and 2012  (edited by M.B. Nathanson), Springer Proc. in Math. Stat.} Vol. 101, Springer, New York, 2014, pp. 257--312.

\bibitem{Sun2020}
Z.-W. Sun.
\newblock New series for powers of $\pi$ and related congruences.
\newblock {\em Electron. Res. Arch.}, 28(3)(2020), 1273--1342.
\bibitem{Sun2020b}
Z.-W. Sun.
\newblock List of conjectural series for powers of $\pi$ and other constants.
\newblock In:{\em  Ramanujan's Identities}, Press of Harbin Institute of Technology, 2020, pp. 205–260.
\bibitem{Sun2021}
Z.-W. Sun.
\newblock New type series for powers of $\pi$.
\newblock {\em J. Comb. Number Theory}, 12(3)(2020), 157--208.\bibitem{Sun2022}
Z.-W. Sun.
\newblock Series with summands involving harmonic numbers.
\newblock {\em M. B. Nathanson (ed.), Combinatorial and Additive Number Theory}, Springer, to appear.

\bibitem{Hamme1997}
L. Van Hamme.
\newblock Some conjectures concerning partial sums of generalized hypergeometric
series.
\newblock in:{\em p-Adic functional analysis (Nijmegen, 1996), Lecture Notes in Pure and Appl. Math.}, 192, Dekker, New York (1997), 223–236.

\bibitem{Hoeven2018}
J. van der Hoeven.
\newblock Creative telescoping using reductions.
\newblock {\em Preprint:hal-01773137v2}, June 2018.


\bibitem{WanZudilin2012}
J. Wan and W. Zudilin.
\newblock Generating functions of Legendre polynomials: A tribute to Fred Brafman.
\newblock {\em J. Approx. Theory}, 164(2012), 488-503.


\bibitem{WZ2023AAM}
R.-H. Wang and M.X.X. Zhong.
\newblock Polynomial reduction for holonomic sequences and applications in $\pi$-series and congruences.
\newblock {\em Adv. in Appl. Math.}, 150(2023), 102568.


\bibitem {Wei2023}
C. Wei.
\newblock On two double series for $\pi$ and their $q$-analogues.
\newblock{\em Ramanujan J.}, 60 (2023),  615--625.


\bibitem{WZ1990}
H.S. Wilf and D. Zeilberger.
\newblock An algorithmic proof theory for hypergeometric (ordinary and ``$q$'') multisum/integral identities.
\newblock {\em Invent. Math.}, 108(1992), 575--633.

\bibitem{Zeilberger1990c}
D. Zeilberger.
\newblock A fast algorithm for proving terminating hypergeometric identities.
\newblock {\em Discrete Math.}, 80(1990), 207--211.


\end{thebibliography}
\end{document}